\documentclass[12pt]{article}

\usepackage{amssymb}
\usepackage{amsmath,amsthm}
\usepackage[dvips]{graphicx}
\usepackage{wrapfig}

\newtheorem{thm}{Theorem}[section]
\newtheorem{cor}[thm]{Corollary}
\newtheorem{lem}[thm]{Lemma}
\newtheorem{prop}[thm]{Proposition}
\newtheorem{q}[thm]{Question}
\theoremstyle{definition}
\newtheorem{defn}[thm]{Definition}
\newtheorem{eg}[thm]{Example}
\theoremstyle{remark}
\newtheorem{rem}[thm]{Remark.}

\title{The warping polynomial of a knot diagram}
\author{Ayaka Shimizu\thanks{Osaka City University Advanced Mathematical Institute, 3-3-138, Sugimoto, Sumiyoshi-ku Osaka 558-8585, Japan. shimizu1984@gmail.com}}
\date{\today}

\begin{document}

\maketitle

\begin{abstract}
We introduce the warping polynomial of an oriented knot diagram. 
In this paper, we characterize the warping polynomial, and define the span of a knot to be the minimal span of the warping polynomial for all diagrams of the knot. 
We show that the span of a knot is one if and only if it is non-trivial and alternating, and we give an inequality between the span and the dealternating number.

\end{abstract}

\section{Introduction}

\noindent Throughout this paper knots and knot diagrams are oriented and diagrams are on $S^2$. 
Kawauchi defined the warping degree which represents a complexity of a knot diagram in \cite{kawauchi}, 
and expanded the concept to spatial graphs in \cite{kawauchi2}. 
The relation between the warping degree and the crossing number of a knot or link diagram is studied by the author in \cite{shimizu} and \cite{shimizu-link}, and recently the warping degree has extended to nanowords by Fukunaga in \cite{fukunaga}. 
In this paper, we define a new Laurent polynomial in one-variable $t$ with non-negative integral coefficients -- the warping polynomial $W_D(t)$ of a knot diagram $D$. 
(We define the warping degree and the warping polynomial in Section \ref{w-poly-sec}.) 
The lower degree of $W_D(t)$ implies the warping degree of $D$, and $W_D(1)/2$ represents the crossing number of $D$. 
Moreover, we have $W_D(-1)=0$ for any knot diagram $D$, and $W_D(0)\ne 0$ if and only if $D$ is a monotone diagram (see Proposition \ref{propwp} and Lemma \ref{minusone}). 
The warping polynomial depends on the orientation: 
Let $-D$ be the diagram $D$ with orientation reversed. 
We have $W_{-D}(t)=t^cW_D(t^{-1})$, where $c$ is the crossing number of $D$ (see Proposition \ref{ori-mirror}). 
In the following theorem, we characterize the warping polynomial:

\phantom{x}
\begin{thm}
A polynomial $f(t)$ is the warping polynomial of a knot diagram if and only if 
\begin{align*}
f(t)=&m_0t^k+(m_0+m_1)t^{k+1}+(m_1+m_2)t^{k+2}+\dots +(m_{l-2}+m_{l-1})t^{k+l-1}+m_{l-1}t^{k+l},
\end{align*}
where $k, l=0,1,2,\dots , m_i=1,2,\dots $ $(i=0,1,2,\dots ,l-1)$ and $m_0+m_1+\dots +m_{l-1}\ge k+l$. 
\label{poly-wd}
\end{thm}
\phantom{x}

\noindent The proof is given in Section \ref{proof-poly-wd}. 
In Section \ref{span-sec}, we define the \textit{span} of a knot diagram $D$, denoted by $\mathrm{spn}(D)$, to be the span of the warping polynomial $W_D(t)$, and also define the \textit{span} of a knot $K$, denoted by $\mathrm{spn}(K)$, to be the minimal $\mathrm{spn}(D)$ for all diagrams $D$ of $K$. 
We have $\mathrm{spn}(K)=0$ if and only if $K$ is the trivial knot, and $\mathrm{spn}(K)=1$ if and only if $K$ is a non-trivial alternating knot (see Theorem \ref{knotspan01}). 
The \textit{dealternating number} $\mathrm{dalt}(D)$ of a knot diagram $D$ is the minimal number of crossing changes which are needed to obtain an alternating diagram from $D$, and the \textit{dealternating number} $\mathrm{dalt}(K)$ of a knot $K$ is the minimal $\mathrm{dalt}(D)$ for all diagrams $D$ of $K$ \cite{ABBCFHJP}. 
We have the following theorem: 

\phantom{x}
\begin{thm}
Let $K$ be a knot. We have 
$$\frac{\mathrm{spn}(K)-1}{2}\le \mathrm{dalt}(K).$$
\label{span-dalt-k}
\end{thm}
\phantom{x}

\noindent The proof is given in Section \ref{span-sec}. 
The rest of this paper is organized as follows: 
In Section \ref{w-poly-sec}, we define the warping polynomial of a knot diagram and show properties. 
In Section \ref{proof-poly-wd}, we study the warping polynomial of a one-bridge diagram and prove Theorem \ref{poly-wd}. 
In Section \ref{span-sec}, we consider the relation between the span and the dealternating number of a knot to prove Theorem \ref{span-dalt-k}, and discuss the spans of connected sums.

\section{Warping polynomial}
\label{w-poly-sec}
\subsection{Warping degree}

\begin{figure}[!ht]
\begin{center}
\includegraphics[width=27mm]{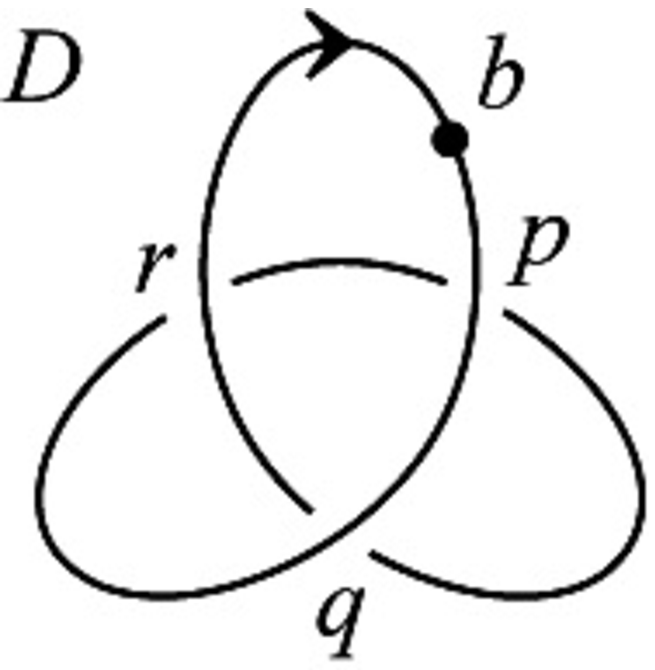}
\caption{}
\label{wcp-ex}
\end{center}
\end{figure}

\begin{figure}[!ht]
\begin{center}
\includegraphics[width=100mm]{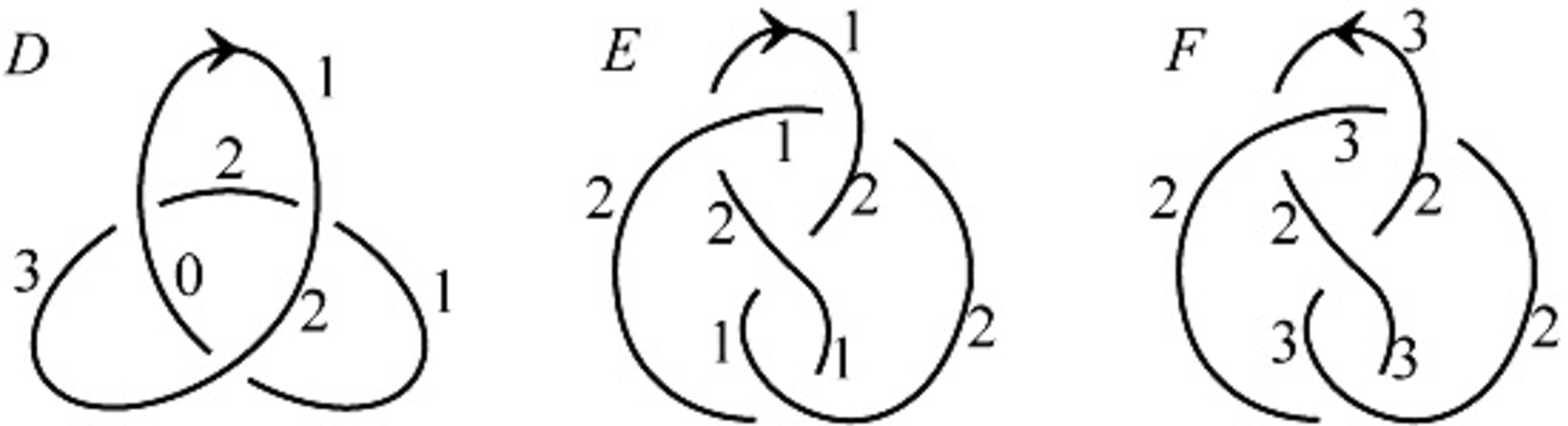}
\caption{}
\label{w-d-l}
\end{center}
\end{figure}

\noindent In this subsection, we explain the warping degree and define the warping degree labeling of a knot diagram. 
A \textit{base point} $b$ of a knot diagram $D$ is a point on $D$ which is not a crossing point. 
We denote the pair of $D$ and $b$ by $D_b$. 
A crossing point of $D$ is a \textit{warping crossing point} of $D_b$ if we come to the point as an under-crossing first 
when we go along $D$ with the orientation by starting from $b$.
The \textit{warping degree} $d(D_b)$ of $D_b$ is the number of warping crossing points of $D_b$ \cite{kawauchi}. 
For example in Figure \ref{wcp-ex}, we have $d(D_b)=1$ because only $r$ is a warping crossing point of $D_b$. 
We note that the similar notions are studied by Fujimura \cite{fujimura}, Fung \cite{fung}, Lickorish and Millett \cite{LM}, Okuda \cite{okuda} and Ozawa \cite{ozawa} considering the ascending number with an orientation. 
An \textit{edge} of $D$ is a path on $D$ between crossing points which has no crossings in the interior. 
We say that an edge $e$ has the warping degree $d$ if $d(D_b)=d$ for a base point $b$ on $e$. 
The \textit{warping degree labeling} of a knot diagram is a labeling of edges with the warping degrees. 
For example, the three diagrams $D$, $E$ and $F$ in Figure \ref{w-d-l} are labeled with warping degree labeling. 
The following lemma helps us to apply the warping degree labeling to a knot diagram: 

\phantom{x}
\begin{lem}{$($Lemma 2.5 in \cite{shimizu}$)$}
For two base points $a_1,a_2$ $($resp. $b_1,b_2)$ which are put across an over-crossing  $($resp. under-crossing$)$ as shown in Figure \ref{lem25}, 
we have $d(D_{a_2})=d(D_{a_1})+1$ $($resp. $d(D_{b_2})=d(D_{b_1})-1)$. 
\begin{figure}[!ht]
\begin{center}
\includegraphics[width=70mm]{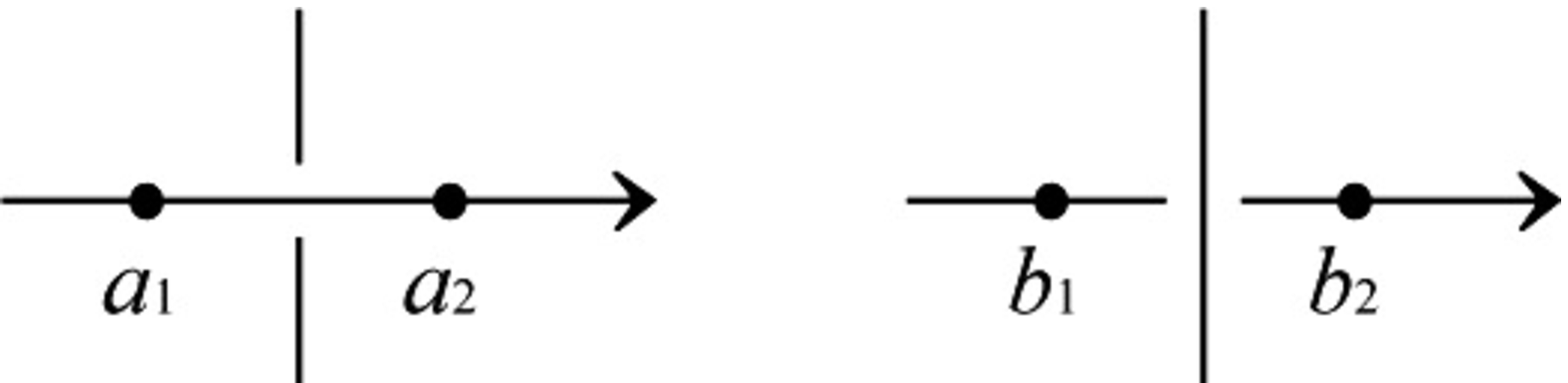}
\caption{}
\label{lem25}
\end{center}
\end{figure}
\label{lemma25}
\end{lem}
\phantom{x}

\noindent The \textit{warping degree} $d(D)$ of a knot diagram $D$ is the minimal $d(D_b)$ for all base points $b$ of $D$ \cite{kawauchi}. 
Let $c(D)$ be the crossing number of $D$. 
In \cite{shimizu}, the following theorem is shown by the author: 

\phantom{x}
\begin{thm}{$($Theorem 1.3 in \cite{shimizu}$)$}
Let $D$ be a knot diagram and let $-D$ be the diagram $D$ with orientation reversed. 
Then we have 
$$d(D)+d(-D)+1\le c(D).$$
Further, the equality holds if and only if $D$ is an alternating diagram. 
\label{knot-wd-thm}
\end{thm}
\phantom{x}

\subsection{Warping polynomial}

\noindent We define the warping polynomial:

\phantom{x}
\begin{defn}
The \textit{warping polynomial} $W_D(t)$ of a knot diagram $D$ is 
$$W_D(t)=\sum _e t^{i(e)},$$
where $i(e)$ is the warping degree of an edge $e$, and $\sum _e$ is the sum for all edges $e$ of $D$. 
\end{defn}
\phantom{x}

\begin{rem}
In other words, the warping polynomial of a knot diagram $D$ is 
$$W_D(t)=\sum _{i\in \mathbb{N}} n_i t^i,$$
where $n_i$ is the number of edges which are labeled $i$ with warping degree labeling.
\end{rem}
\phantom{x}

\noindent For example, the diagrams in Figure \ref{w-d-l} have the warping polynomials $W_D(t)=1+2t+2t^2+t^3$, $W_E(t)=4t+4t^2$ and $W_F(t)=4t^2+4t^3$. 
We have the following proposition:

\phantom{x}
\begin{prop}
Let $D$ be an alternating knot diagram with $c(D)\ge 1$. 
Then, $W_D(t)=ct^d+ct^{d+1}$, where $c=c(D)$ and $d=d(D)$. 

\phantom{x}
\begin{proof}
When we go along $D$, we go through edges labeled $d$ and $d+1$ alternately. 
Since $D$ has $2c$ edges, we have $W_D(t)=ct^d+ct^{d+1}$. 
\end{proof}
\end{prop}
\phantom{x}

\noindent A knot diagram $D$ is \textit{monotone} if $d(D)=0$ \cite{kawauchi}. 
We remark that a monotone diagram represents the trivial knot. 
Let $\mathrm{l}\text{-}\mathrm{deg} W_D(t)$ (resp. $\mathrm{u}\text{-}\mathrm{deg} W_D(t)$) be the lower (resp. upper) degree of $W_D(t)$. 
For example, $\mathrm{l}\text{-}\mathrm{deg} W_D(t)=0$ and $\mathrm{u}\text{-}\mathrm{deg} W_D(t)=3$ for $W_D(t)=1+2t+2t^2+t^3$. 
We have the following proposition:

\phantom{x}
\begin{prop}
The warping polynomial has the following properties:
\begin{description}
 \item[(i)] $W_{\bigcirc}(t)=1$, where $\bigcirc$ is the knot diagram with $c(\bigcirc )=0$, 
 \item[(ii)] $\mathrm{l}\text{-}\mathrm{deg} W_D(t)=d(D)$, 
 \item[(iii)] $D$ is monotone if and only if $W_D(0)\ne 0$, and 
 \item[(iv)] if $c(D)\ge 1$, then $W_D(1)=2c(D)$.
\end{description}

\phantom{x}
\begin{proof}
\begin{description}
 \item[(i)--(iii)] Obvious by definition. 
 \item[(iv)] By definition, $W_D(1)=\sum _e 1=\sharp \{ \text{edges of } D \} =2c(D)$. 
\end{description}
\end{proof}
\label{propwp}
\end{prop}
\phantom{x}

\noindent From Lemma \ref{lemma25}, we have the following proposition: 

\phantom{x}
\begin{prop}
Let $W_D(t)=n_kt^k+n_{k+1}t^{k+1}+\dots +n_{k+l}t^{k+l}$ be the warping polynomial of a knot diagram $D$ with $n_k,n_{k+l}\ne 0$ $(0\le k, l \le c(D))$. 
Then, $n_{k+1},n_{k+2},\dots ,n_{k+l-1}\ne 0$. 
\label{by25}
\end{prop}
\phantom{x}

\noindent For $-D$ and the mirror image $D^*$ of a knot diagram $D$, we have the following proposition: 

\phantom{x}
\begin{prop}
We have 
$$W_{-D}(t)=W_{D^*}(t)=t^cW_D(t^{-1}), $$
where $c=c(D)$. 

\phantom{x}
\begin{proof}
Let $b$ be a base point of $D$. 
A crossing point $p$ is a warping crossing point of $-D_b$ if and only if $p$ is not a warping crossing point of $D_b$. 
Hence $d(D_b)+d(-D_b)=c(D)$ (Lemma 2.1 in \cite{shimizu}). 
Then we have 
$$W_{-D}(t)=\sum _e t^{j(e)}=\sum _e t^{c-i(e)}=t^c W_D(t^{-1}),$$
where $i(e)$ represents the warping degree labeling of $D$ and $j(e)$ represents that of $-D$. 
Similarly, we have $d(D_b)+d(D_b^*)=c(D)$ for the mirror image $D^*$ of $D$, and we have $W_{D^*}(t)=t^cW_D(t^{-1})$. 
\end{proof}
\label{ori-mirror}
\end{prop}
\phantom{x}

\noindent We have the following lemma:

\phantom{x}
\begin{lem}
For a knot diagram $D$ with $c(D)\ge 1$, we have 
$$W_D(-1)=0.$$

\phantom{x}
\begin{proof}
From Lemma \ref{lemma25}, when we go along $D$ with warping degree labeling, we go through edges with even warping degrees and odd warping degrees alternately. 
Hence the number of edges labeled with an odd number is equal to that with an even number. 
Then, for the warping polynomial $W_D(t)=\sum _{i\in \mathbb{N}} n_i t^i$, we have 
$n_d+n_{d+1}+\dots =n_{d+1}+n_{d+3}+\dots $, and therefore $n_d-n_{d+1}+n_{d+2}-n_{d+3}+-\dots =(-1)^d W_D(-1)=0$. 
\end{proof}
\label{minusone}
\end{lem}
\phantom{x}

\noindent We have the following corollary: 

\phantom{x}
\begin{cor}
Let $W_D(t)=\sum _{i\in \mathbb{N}} n_i t^i$ be a warping polynomial. Then, 
$$\sum _{i: \text{ odd}} n_i=\sum _{i: \text{ even}} n_i=c(D).$$

\phantom{x}
\begin{proof}
From Lemma \ref{minusone}, we have the first equality. 
From Proposition \ref{propwp} (iv), $\sum _i n_i=2c(D)$. 
\end{proof}
\end{cor}
\phantom{x}

\subsection{Reidemeister moves}

\noindent In this subsection, we consider the warping polynomials on Reidemaister moves. 
Polyak showed in \cite{polyak} the following theorem:

\phantom{x}
\begin{thm}{$($Theorem 1.1 in \cite{polyak}$)$}
Let $D$ and $D'$ be two diagrams in $\mathbb{R}^2$, representing the same oriented link. 
Then one may pass from $D$ to $D'$ by isotopy and a finite sequence of four oriented Reidemeister moves $\Omega _{1a}$, $\Omega _{1b}$, $\Omega _{2a}$, and $\Omega _{3a}$, shown in Figure \ref{ori-rei}. 
\begin{figure}[!ht]
\begin{center}
\includegraphics[width=120mm]{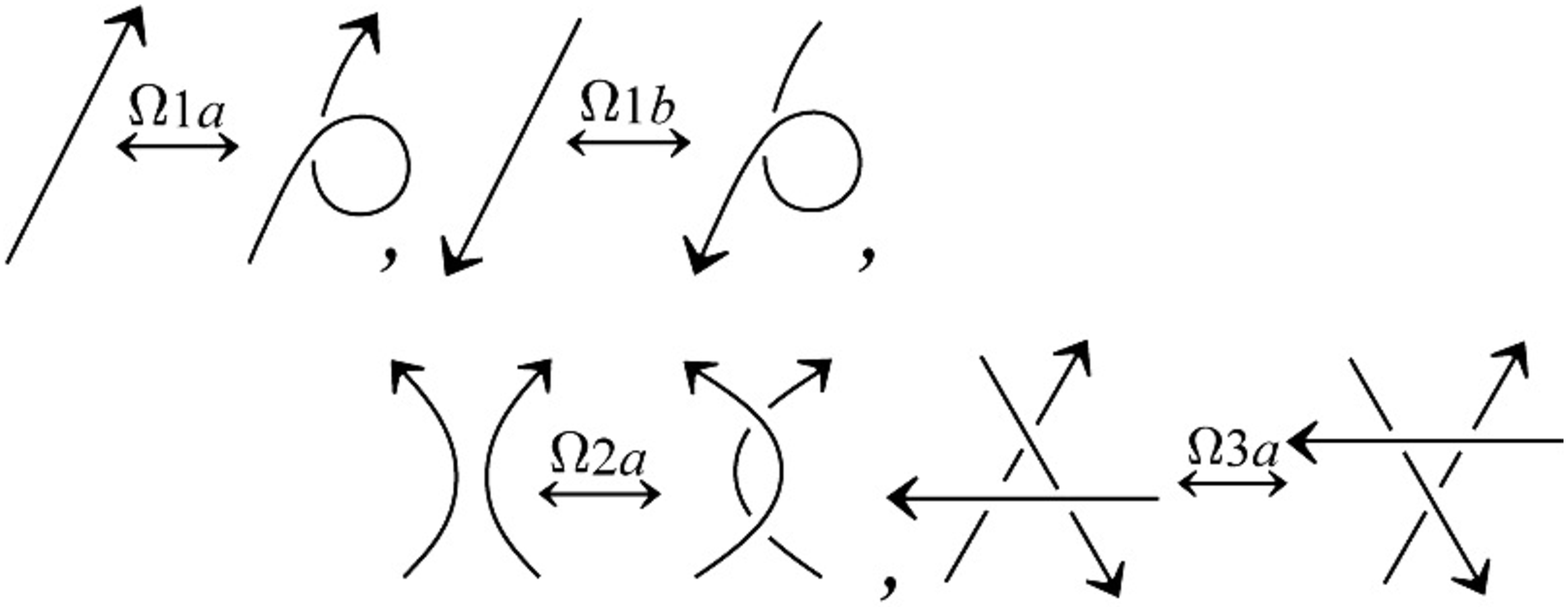}
\caption{}
\label{ori-rei}
\end{center}
\end{figure}
\end{thm}
\phantom{x}

\begin{figure}[!ht]
\begin{center}
\includegraphics[width=120mm]{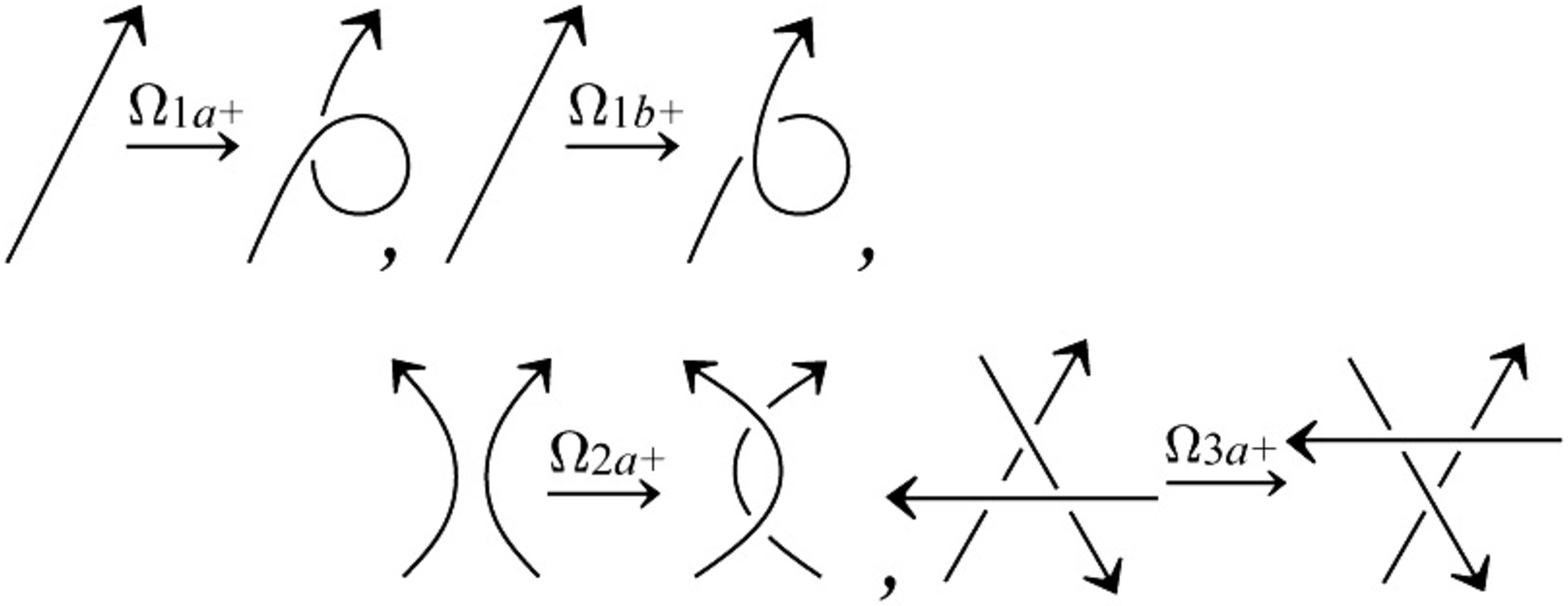}
\caption{}
\label{ori-rei-ori}
\end{center}
\end{figure}

\noindent In this paper, we call the four oriented Reidemeister moves in Figure \ref{ori-rei-ori} 
$\Omega _{1a+}$, $\Omega _{1b+}$, $\Omega _{2a+}$, and $\Omega _{3a+}$. 
With respect to the Reidemeister moves of type $\Omega _{1a+}$ and $\Omega _{1b+}$, we have the following lemma: 

\phantom{x}
\begin{lem}{$($Proposition 8 in \cite{shimizu-j}$)$}
Let $D$ be a knot diagram. 
Let $D'$ (resp. $D''$) be a knot diagram which is obtained from $D$ by a Reidemeister move of type $\Omega _{1a+}$ (resp. $\Omega _{1b+}$) 
at an edge of $D$ whose warping degree is $i$ (Figure \ref{r-1-fig}). 
Then we have 
\begin{align*}
W_{D'}(t)&=W_D(t)+t^i(1+t),\\
W_{D''}(t)&=tW_D(t)+t^i(1+t).
\end{align*}
\begin{figure}[!ht]
\begin{center}
\includegraphics[width=90mm]{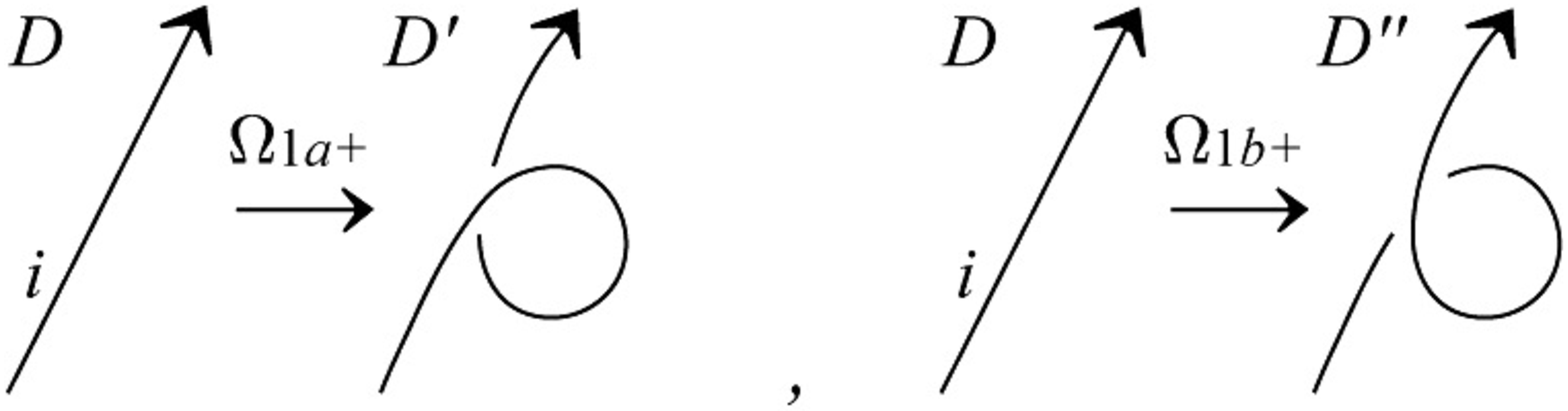}
\caption{}
\label{r-1-fig}
\end{center}
\end{figure}

\phantom{x}
\begin{proof}
Let $a, a_i, b, b_i$ $(i=1,2,3)$ be the base points and $p, q$ the crossing points as shown in Figure \ref{r-1-pf}. 
\begin{figure}[!ht]
\begin{center}
\includegraphics[width=90mm]{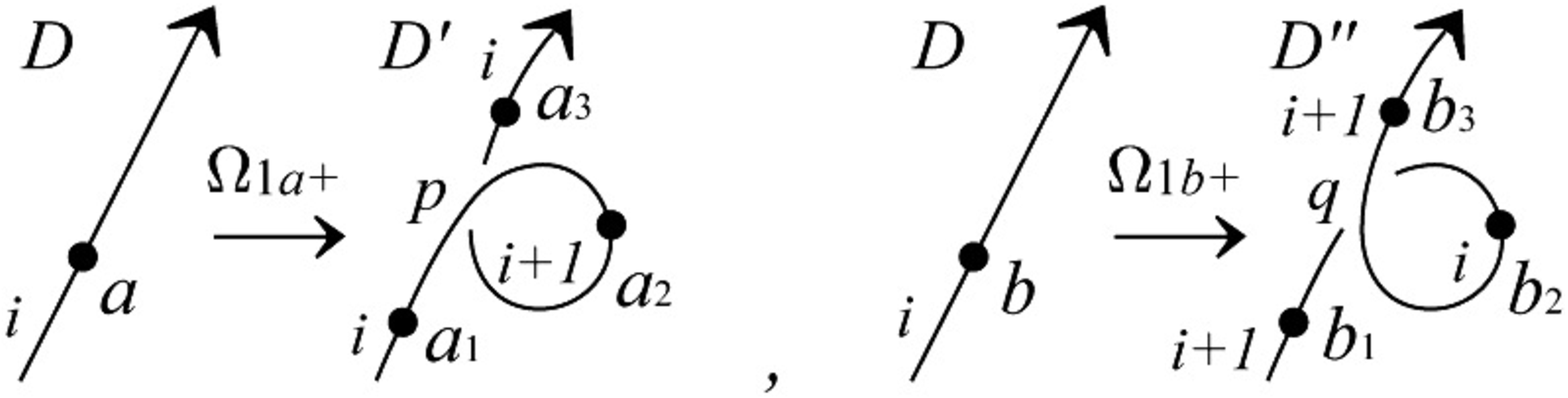}
\caption{}
\label{r-1-pf}
\end{center}
\end{figure}
Then, $D'_{a_1}$ and $D'_{a_3}$ have warping crossing points at the same crossing points as $D_a$. 
Since $p$ is a warping crossing point of $D'_{a_2}$, $D'_{a_2}$ has one more warping crossing point than $D'_{a_1}$, $D'_{a_3}$, and therefore $D_a$ as well. 
Then, $W_{D'}(t)=W_D(t)+t^i+t^{i+1}$. 
Similarly, $D''_{b_1}$ and $D''_{b_3}$ have one more warping crossing point than $D''_{b_2}$, 
and $D''_{b_2}$ has warping crossing points at the same crossing points as $D_b$. 
Then, $W_{D''}(t)=tW_D(t)+t^i+t^{i+1}$. 
\end{proof}
\label{r-1-lem}
\end{lem}
\phantom{x}

\noindent We show an example:

\phantom{x}
\begin{eg}
For the knot diagrams in Figure \ref{r-1-ex-fig}, we have 
\begin{align*}
W_D(t)&=3t+3t^2,\\
W_E(t)&=4t+4t^2=W_D(t)+t(1+t),\\
W_F(t)&=3t+4t^2+t^3=W_D(t)+t^2(1+t),\\
W_G(t)&=t+4t^2+3t^3=tW_D(t)+t(1+t).\\
\end{align*}
\begin{figure}[!ht]
\begin{center}
\includegraphics[width=120mm]{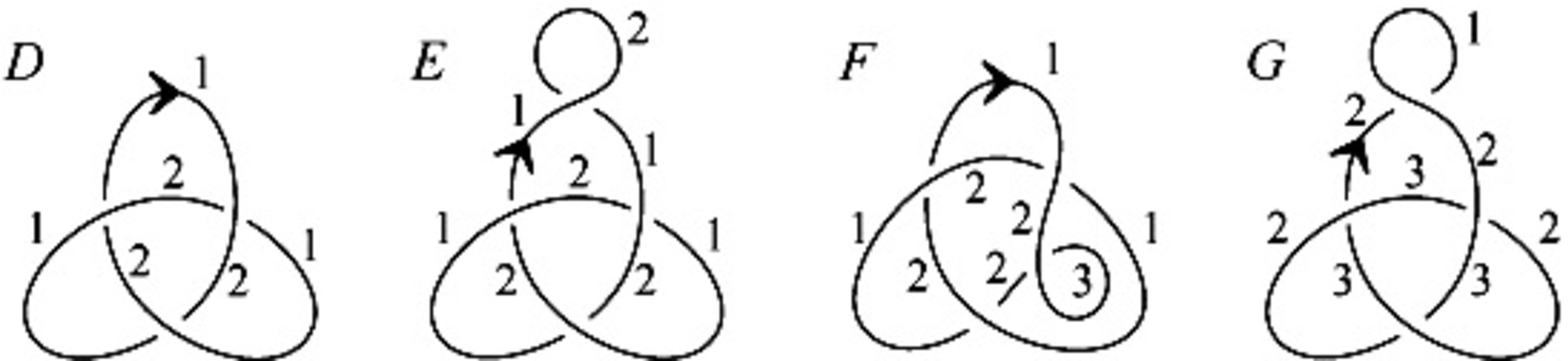}
\caption{}
\label{r-1-ex-fig}
\end{center}
\end{figure}
\end{eg}
\phantom{x}

\noindent With respect to the Ridemeister moves of type $\Omega _{2a}$ and $\Omega _{3a}$, see Propositions 10 and 11 in \cite{shimizu-j}.

\section{Proof of Theorem \ref{poly-wd}}
\label{proof-poly-wd}

\begin{figure}[!ht]
\begin{center}
\includegraphics[width=55mm]{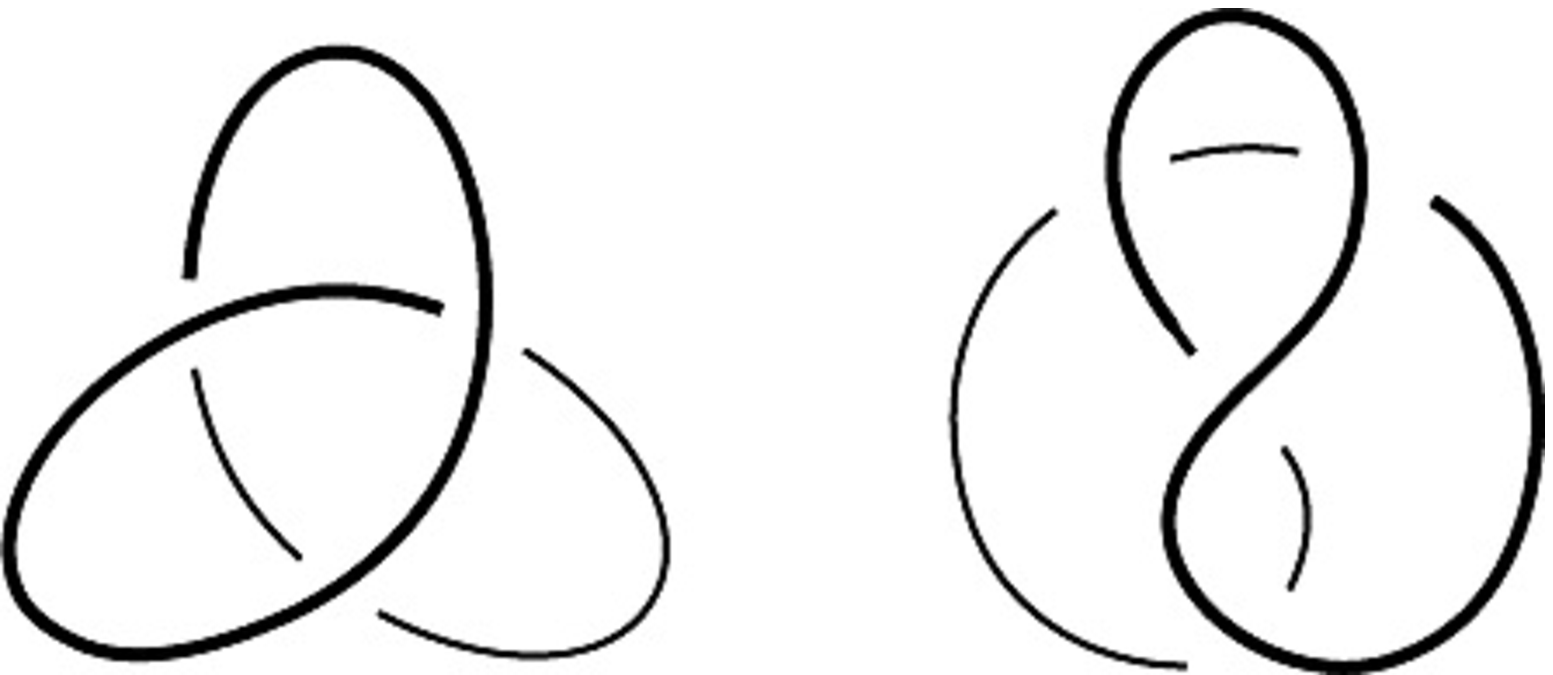}
\caption{}
\label{onebridge}
\end{center}
\end{figure}

\noindent In this section, we prove Theorem \ref{poly-wd}. 
A \textit{bridge} in a knot diagram $D$ is a path on $D$ between under-crossings which has no under-crossings and at least one over-crossing in the interior. 
A knot diagram $D$ is a \textit{one-bridge diagram} if $D$ has exactly one bridge (Figure \ref{onebridge}). 
Let $F^l(t)=1+2t+2t^2+\dots +2t^{l-1}+t^l$. 
We have the following lemma:

\phantom{x}
\begin{lem}
A knot diagram $D$ with $c(D)=k$ $(k=1,2,\dots )$ is a one-bridge diagram if and only if $W_D(t)=F^k(t).$

\phantom{x}
\begin{proof}
Let $D$ be a knot diagram with $c(D)=k$. 
By Lemma \ref{lemma25}, the following four conditions are equivalent: 
\begin{itemize}
\item $D$ is a one-bridge diagram. 
\item When going along $D$, we come to $k$ over-crossings in a row and then $k$ under-crossings in a row. 
\item $D$ has just one edge with warping degree $0$, two edges with warping degree $1, 2, \dots ,k-1$, respectively, and one edge with warping degree $k$. 
\item $W_D(t)=1+2t+2t^2+\dots +2t^{k-1}+t^k$.
\end{itemize}
\end{proof}
\label{one-b-lem}
\end{lem}
\phantom{x}

\noindent We prove Theorem \ref{poly-wd}. \\

\begin{figure}[!ht]
\begin{center}
\includegraphics[width=70mm]{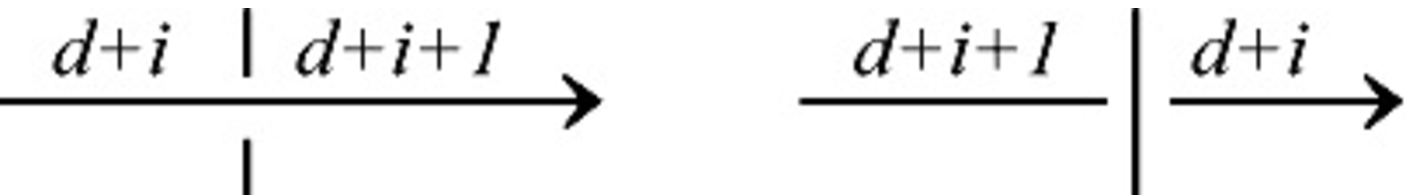}
\caption{}
\label{pf11}
\end{center}
\end{figure}

\phantom{x}
\noindent Proof of Theorem \ref{poly-wd}. 
Let $D$ be a knot diagram with $d(D)=d$ and $\max _b d(D_b)=d+s$ ($d, s=0,1,2,\dots $). 
Let $n_i$ be the number of crossing points of $D$ with warping degree labeling as depicted in the left-hand side of Figure \ref{pf11} ($i=0, 1, 2, \dots , s-1, n_i\ge 1$). 
Then $D$ also has $n_i$ crossing points as shown in the right-hand side of Figure \ref{pf11}. 
Therefore, the number of edges of $D$ which has the warping degree $d+j$ is $n_{j-1}+n_j$ because the number of edges with warping degree $d+j$ whose start-points are over-crossings (resp. under-crossings) is $n_{j-1}$ (resp. $n_j$). 
Hence we have $W_D(t)=n_0t^d+(n_0+n_1)t^{d+1}+\dots +(n_{s-2}+n_{s-1})t^{d+s-1}+n_{s-1}t^{d+s}$. 
Since $W_D(1)=2(n_0+n_1+\dots +n_{s-1})=2c(D)$ and $\max _b d(D_b)=d+s\le c(D)$, we have $n_0+n_1+\dots +n_{l-1}\ge d+s$. \\
Let $f(t)=m_0t^k+(m_0+m_1)t^{k+1}+(m_1+m_2)t^{k+2}+\dots +(m_{l-2}+m_{l-1})t^{k+l-1}+m_{l-1}t^{k+l}$ be a polynomial 
with $k, l=0,1,2,\dots , m_i=1,2,\dots $ $(i=0,1,2,\dots ,l-1)$ and $m_0+m_1+\dots +m_{l-1}\ge k+l$. 
Then we can divide each $m_i$ into ${m_i}'$ and ${m_i}''$ so that $m_i={m_i}'+{m_i}''+1$ and ${m_0}'+{m_1}'+\dots +{m_{l-1}}'=k$. 
Let $E$ be a one-bridge knot diagram with $c(E)=l$, i.e., $W_E(t)=F^l(t)$. 
We apply $\Omega _{1b+}$ to $E$ at the edges whose warping degrees are $i$ ${m_i}'$ times for all $i$. 
Then we obtain a knot diagram $E'$ with $W_{E'}(t)=t^k F^l (t)+\sum _{i=0}^{l-1} {m_i}'t^{k+i}(1+t)$. 
We apply $\Omega _{1a+}$ to $E'$ at the edges whose warping degrees are $k+i$ ${m_i}''$ times for all $i$. 
Thus, we obtain a knot diagram $E''$ with 
\begin{align*}
W_{E''}(t)&=t^k F^l (t)+\sum _{i=0}^{l-1} {m_i}'t^{k+i}(1+t)+\sum _{i=0}^{l-1} {m_i}''t^{k+i}(1+t)\\
            &=f(t).
\end{align*}
\hfill$\square$
\phantom{x}

\section{The span}
\label{span-sec}

\noindent Let $\mathrm{span}f(t)$ be the span of a polynomial $f(t)$, namely $\mathrm{span}f(t)=\mathrm{u}\text{-}\mathrm{deg} f(t)-\mathrm{l}\text{-}\mathrm{deg} f(t)$. 
As we have seen, an alternating knot diagram $D$ with $c(D)\ge 1$ has $\mathrm{span}W_D(t)=1$ and 
a one-bridge diagram $E$ with $c(E)=l$ has $\mathrm{span}W_E(t)=l$. 
In this section, we study the spans of a knot diagram and a knot, and prove Theorem \ref{span-dalt-k}. 
We also discuss the spans of connected sums. 

\subsection{The span and the dealternating number}

\begin{figure}[!ht]
\begin{center}
\includegraphics[width=100mm]{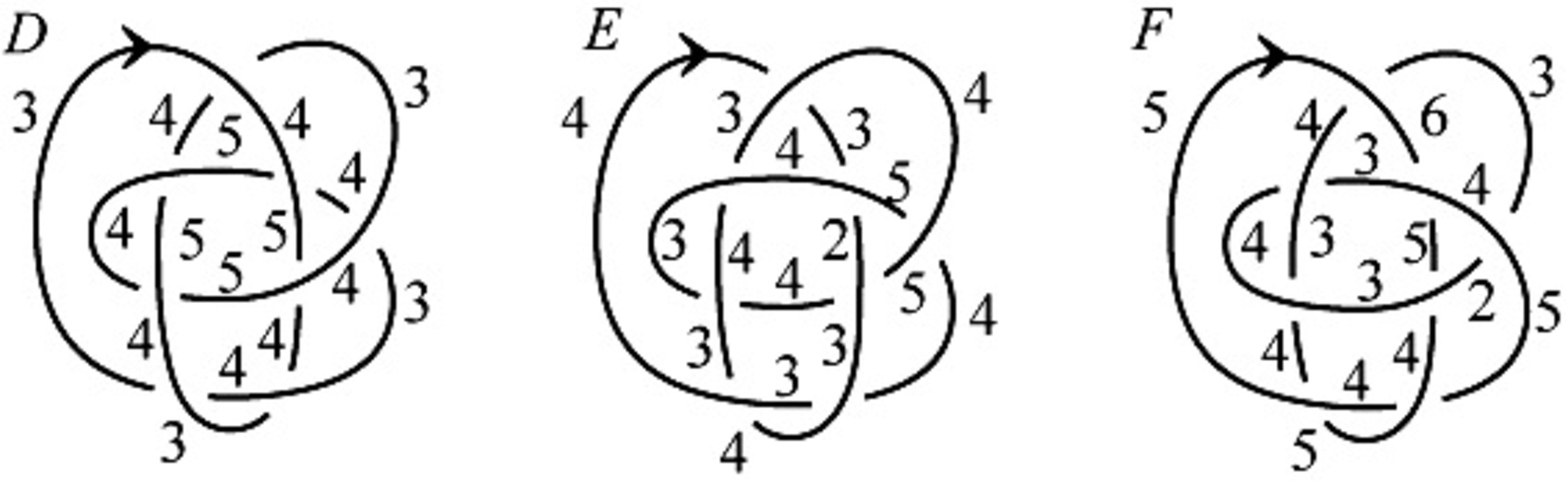}
\caption{}
\label{span8}
\end{center}
\end{figure}

\noindent The \textit{span} $\mathrm{spn}(D)$ of a knot diagram $D$ is defined to be $\mathrm{span}W_D(t)$. 
Obviously, we have $\mathrm{spn}(D)=\max _b d(D_b)-\min _b d(D_b)$. 
For example, the knot diagrams $D$, $E$ and $F$ in Figure \ref{span8} which are diagrams of $8_{19}$, $8_{20}$ and $8_{21}$ have $\mathrm{spn}(D)=2$, $\mathrm{spn}(E)=3$ and $\mathrm{spn}(F)=4$. 
We also have $\mathrm{spn}(D)=c(D)-(d(D)+d(-D))$ because $d(D_b)+d(-D_b)=c(D)$ and $\max _b d(D_b)=c(D)-d(-D)$. 
Hence $\mathrm{spn}(D)$ does not depend on the orientation of $D$. 
We note that $c(D)-(d(D)+d(-D))$ is also discussed in \cite{shimizu}. 
From Theorem \ref{knot-wd-thm}, we immediately have the following corollary: 

\phantom{x}
\begin{cor}
Let $D$ be a knot diagram. Then, 
\begin{itemize}
\item $D$ is the diagram with $c(D)=0$ if and only if $\mathrm{spn}(D)=0$, 
\item $D$ is an alternating diagram with $c(D)\ge 1$ if and only if $\mathrm{spn}(D)=1$. 
\end{itemize}
\label{span-one-cor}
\end{cor}
\phantom{x}

\noindent The \textit{span} $\mathrm{spn}(K)$ of a knot $K$ is the minimal $\mathrm{spn}(D)$ for all diagrams $D$ of $K$. 
From Corollary \ref{span-one-cor}, we have the following theorem: 

\phantom{x}
\begin{thm}
Let $K$ be a knot. Then, 
\begin{itemize}
\item $K$ is the trivial knot if and only if $\mathrm{spn}(K)=0$, 
\item $K$ is a non-trivial alternating knot if and only if $\mathrm{spn}(K)=1$. 
\end{itemize}

\phantom{x}
\begin{proof}
A knot $K$ is trivial if and only if $K$ has a diagram $D$ such that $c(D)=0$, namely $\mathrm{spn}(D)=0$ and $\mathrm{spn}(K)=0$. 
Let $L$ be a non-trivial knot. 
Then $L$ is alternating if and only if $L$ has an alternating diagram $E$ with $c(E)\ge 1$, namely $\mathrm{spn}(E)=1$ and $\mathrm{spn}(L)=1$. 
\end{proof}
\label{knotspan01}
\end{thm}
\phantom{x}

\noindent We have the following proposition: 

\phantom{x}
\begin{prop}
Let $K$ be a knot. 
If $K$ has a closed 3-braid diagram with all positive (or negative) crossings, then $\mathrm{spn}(K)\le 2$. 
Moreover, we have $\mathrm{spn}(K)=2$ if $K$ is neither trivial nor alternating. 

\phantom{x}
\begin{proof}
A closed $n$-braid knot diagram with all positive (or negative) crossings has the span $n-1$ (see Figure \ref{p3-braid}). 
\begin{figure}[!ht]
\begin{center}
\includegraphics[width=25mm]{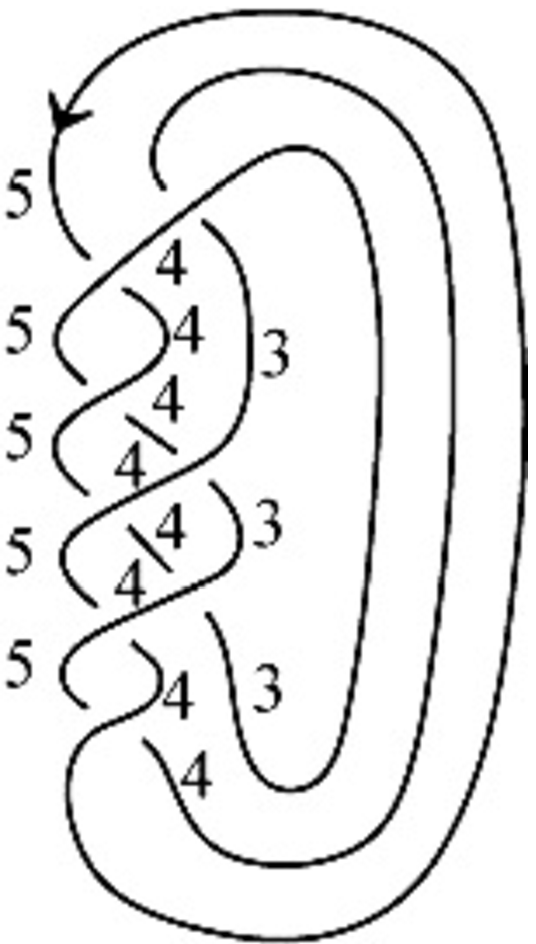}
\caption{}
\label{p3-braid}
\end{center}
\end{figure}
\end{proof}
\end{prop}
\phantom{x}

\noindent With respect to a crossing change, we have the following lemma: 

\phantom{x}
\begin{lem}
Let $D$ be a knot diagram with $c(D)\ge 1$, and let $D'$ be a knot diagram which is obtained from $D$ by a crossing change. 
Then, we have 
$$\vert \mathrm{spn}(D')-\mathrm{spn}(D)\vert \le 2.$$

\phantom{x}
\begin{proof}
Let $p$ be a crossing point of a knot diagram $D$, and let $e_1, e_2, \dots ,e_{2c}$ be the edges of $D$ $(c=c(D))$. 
We divide the edges into two sets $E_1$ and $E_2$ so that they satisfy the following: 
\begin{itemize}
\item If we go through $e_i$ when we go along $D$ from the over-crossing of $p$ to the under-crossing of $p$, then $e_i\in E_1$. 
\item If we go through $e_i$ when we go along $D$ from the under-crossing of $p$ to the over-crossing of $p$, then $e_i\in E_2$. 
\end{itemize}
We define polynomials $f_D(t)$ and $g_D(t)$ as follows: 
$$f_D(t)=\sum_{e_j\in E_1}t^{i(e_j)},  \  g_D(t)=\sum _{e_j\in E_2}t^{i(e_j)},$$
where $i(e_j)$ is the warping degree of $e_j$. 
We have $W_D(t)=f_D(t)+g_D(t)$ by definition. 
Let $D'$ be the knot diagram which is obtained from $D$ by the crossing change at $p$. 
Then we have $f_{D'}(t)=tg_D(t)$, $g_{D'}=t^{-1}f_D(t)$, and therefore $W_{D'}(t)=tg_D(t)+t^{-1}f_D(t)$. 
\end{proof}
\label{cc-2}
\end{lem}

\noindent We show an example: 

\phantom{x}
\begin{eg}
For the knot diagrams $D$ and $D'$ in Figure \ref{f-g-ex}, we obtain $D'$ from $D$ by the crossing change at the crossing point $p$ of $D$. 
\begin{figure}[!ht]
\begin{center}
\includegraphics[width=70mm]{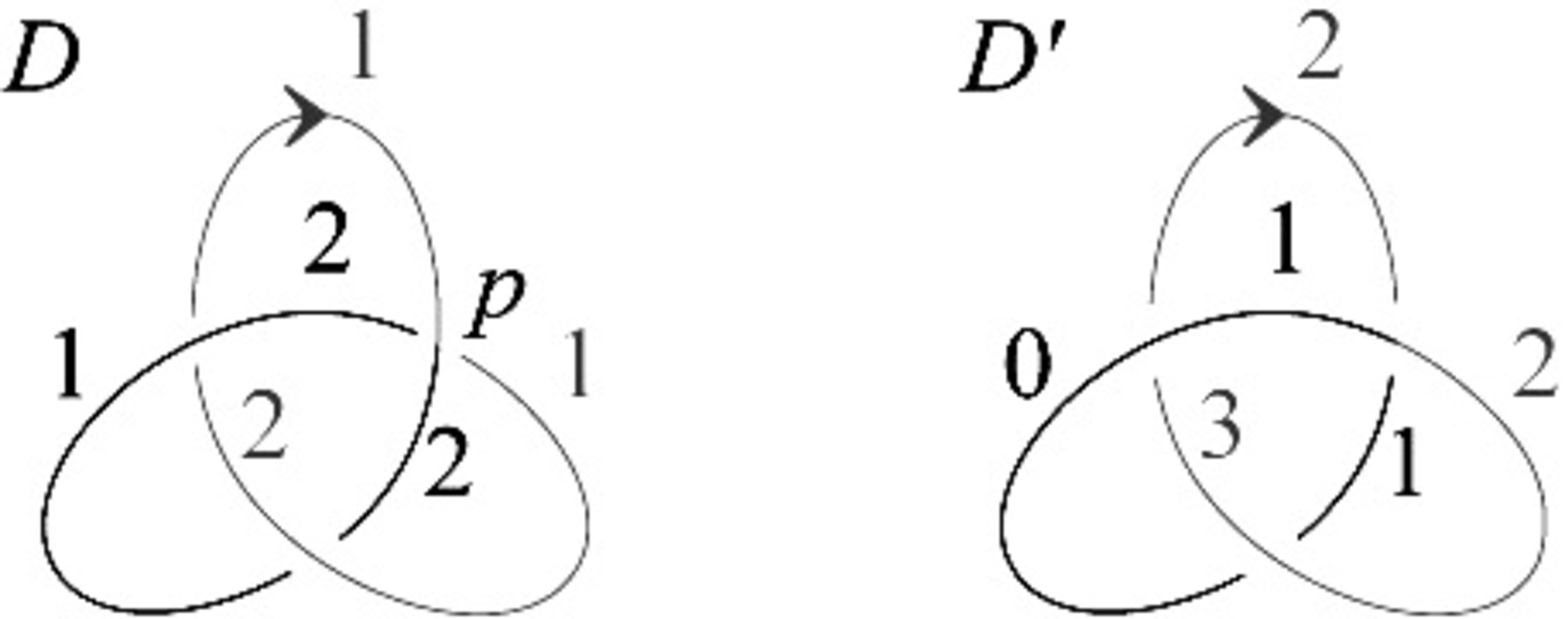}
\caption{}
\label{f-g-ex}
\end{center}
\end{figure}
We have $W_D(t)=3t+3t^2$, $f_D(t)=t+2t^2$ and $g_D(t)=2t+t^2$. 
Then we have $W_{D'}(t)=1+2t+t+2t^2+t^3=t^{-1}f_D(t)+tg_D(t)$, and $\mathrm{spn}(D')-\mathrm{spn}(D)=3-1=2$. 
\end{eg}
\phantom{x}

\noindent We have the following lemma: 

\phantom{x}
\begin{lem}
Let $D$ be a knot diagram. We have 
$$\frac{\mathrm{spn}(D)-1}{2}\le \mathrm{dalt}(D).$$

\phantom{x}
\begin{proof}
Let $\mathrm{dalt}(D)=k$. 
Since the span of an alternating knot diagram with at least one crossing is $1$, we have $\mathrm{spn}(D)\le 1+2k$ by Lemma \ref{cc-2}. 
\end{proof}
\label{span-dalt}
\end{lem}
\phantom{x}

\noindent From the well-known inequality $\mathrm{dalt}(D)\le c(D)/2$ and Lemma \ref{span-dalt}, we can determine the dealternating numbers by $\mathrm{spn}(D)$ and $c(D)$ for some knot diagrams $D$: 

\phantom{x}
\begin{eg}
We have $\mathrm{spn}(D)=8$ and $c(D)=9$ for the diagram $D$ in Figure \ref{9-dalt-ex}, and therefore $7/2\le \mathrm{dalt}(D) \le 9/2$. 
Hence we have $\mathrm{dalt}(D)=4$. 
\begin{figure}[!ht]
\begin{center}
\includegraphics[width=33mm]{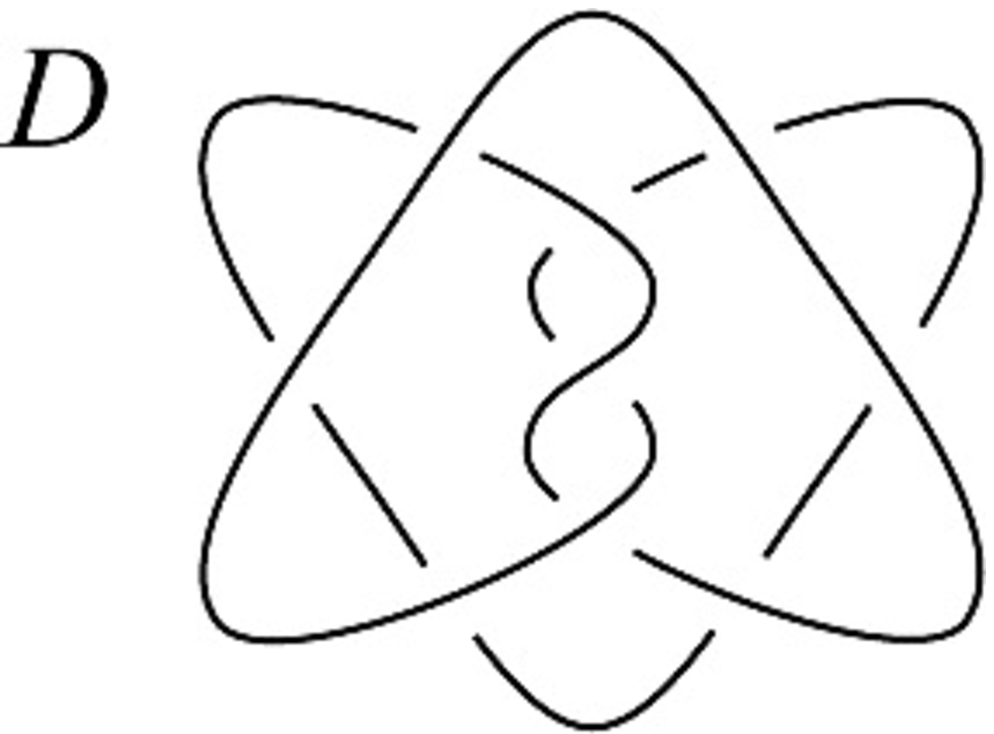}
\caption{}
\label{9-dalt-ex}
\end{center}
\end{figure}
\end{eg}
\phantom{x}

\noindent A knot diagram $D$ is \textit{almost alternating} if $\mathrm{dalt}(D)=1$ \cite{ABBCFHJP}. 
We have the following corollary: 

\phantom{x}
\begin{cor}
Let $D$ be an almost alternating knot diagram. 
Then, $\mathrm{spn}(D)$ is $2$ or $3$. 

\phantom{x}
\begin{proof}
From Theorem \ref{span-dalt}, we have $(\mathrm{spn}(D)-1)/2\le \mathrm{dalt}(D)=1$, i.e., $\mathrm{spn}(D)\le 3$. 
Since $D$ is not alternating, the span is $2$ or $3$. 
\end{proof}
\end{cor}
\phantom{x}

\noindent Further, for almost alternating knot diagrams, we have the following proposition: 

\phantom{x}
\begin{prop}
Let $D$ be an almost alternating knot diagram. 
If $D$ is a knot diagram obtained from an alternating knot diagram by a Reidemeister move of type $\Omega _{1a+}$ or $\Omega _{1b+}$, then $\mathrm{spn}(D)=2$. 
Otherwise, $\mathrm{spn}(D)=3$. 

\phantom{x}
\begin{proof}
Let $D$ be the almost alternating knot diagram obtained from an alternating diagram by the crossing change at a crossing point $p$. 
Let $f_D(t)$ and $g_D(t)$ be polynomials with respect to $p$ as we defined in the proof of Lemma \ref{cc-2}. 
We remark that the span of at least one of $f_D(t)$ and $g_D(t)$ is not zero because $c(D)\ge 2$. 
Since $D$ is alternating except around $p$, we have $\mathrm{span}f_D(t)\le 1$ and $\mathrm{span}g_D(t)\le 1$. 
The edge just after (resp. before) the over-crossing of $p$ has the warping degree $\mathrm{l}\text{-}\mathrm{deg} f_D(t)$ (resp. $\mathrm{u}\text{-}\mathrm{deg} g_D(t)$). 
Hence we have $\mathrm{l}\text{-}\mathrm{deg} f_D(t)=\mathrm{u}\text{-}\mathrm{deg} g_D(t)+1$. 
Therefore, $\mathrm{span}W_D(t)=\mathrm{span}(f_D(t)+g_D(t))=2$ if and only if $\mathrm{span}f_D(t)=0$ or $\mathrm{span}g_D(t)=0$, and if and only if $p$ is the crossing shown in Figure \ref{r1c}. 
\begin{figure}[!ht]
\begin{center}
\includegraphics[width=40mm]{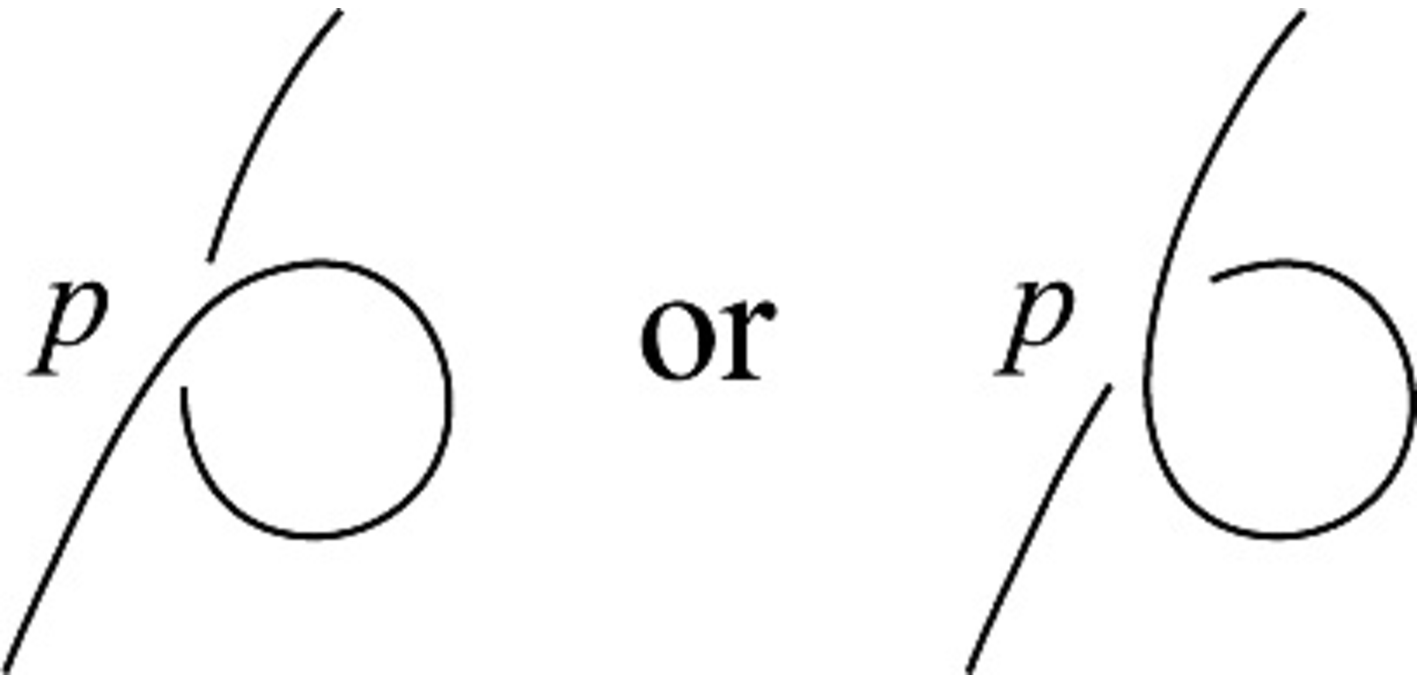}
\caption{}
\label{r1c}
\end{center}
\end{figure}
\end{proof}
\label{aad23}
\end{prop}
\phantom{x}

\noindent We show an example: 

\phantom{x}
\begin{eg}
For two diagrams $D$ and $E$ in Figure \ref{a-alt-ex}, we have $W_D(t)=t+4t^2+3t^3$, $W_E(t)=1+2t+2t^2+t^3$, 
and $\mathrm{spn}(D)=2$, $\mathrm{spn}(E)=3$. 
\begin{figure}[!ht]
\begin{center}
\includegraphics[width=65mm]{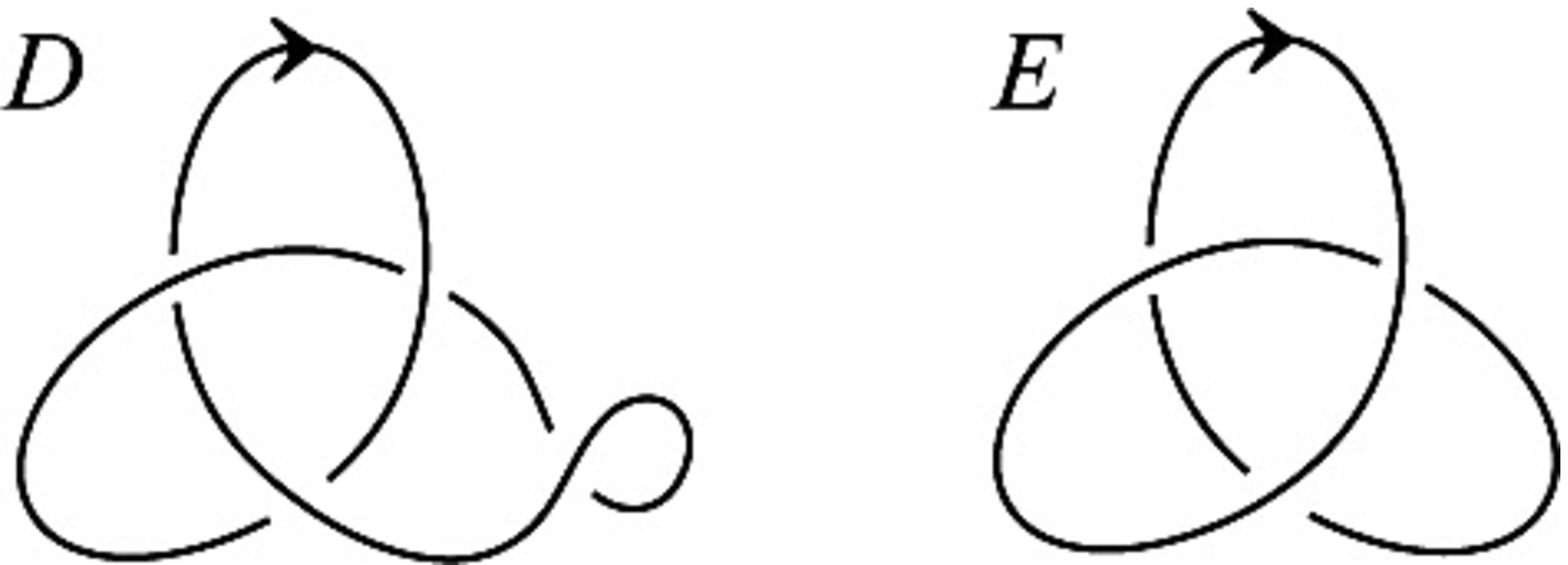}
\caption{}
\label{a-alt-ex}
\end{center}
\end{figure}
\end{eg}
\phantom{x}

\noindent A knot $K$ is \textit{almost alternating} if $K$ has an almost alternating diagram and $K$ is not alternating \cite{ABBCFHJP}. 
We have the following corollary: 

\phantom{x}
\begin{cor}
If $K$ is an almost alternating knot, then $\mathrm{spn}(K)$ is $2$ or $3$. 

\phantom{x}
\begin{proof}
An almost alternating knot $K$ has a diagram $D$ with $\mathrm{spn}(D)=3$ by Proposition \ref{aad23}. 
Since $K$ is neither trivial nor alternating, the span is $2$ or $3$. 
\end{proof}
\end{cor}
\phantom{x}

\noindent Now we prove Theorem \ref{span-dalt-k}. 

\phantom{x}
\noindent Proof of Theorem \ref{span-dalt-k}. 
Let $D$ be a diagram of a knot $K$ with $\mathrm{dalt}(D)=\mathrm{dalt}(K)$. 
Then, $\mathrm{spn}(K)\le \mathrm{spn}(D)\le 2\mathrm{dalt}(D)+1=2\mathrm{dalt}(K)+1$.\\
\hfill$\square$
\phantom{x}

\subsection{Connected sum}

\begin{figure}[!ht]
\begin{center}
\includegraphics[width=130mm]{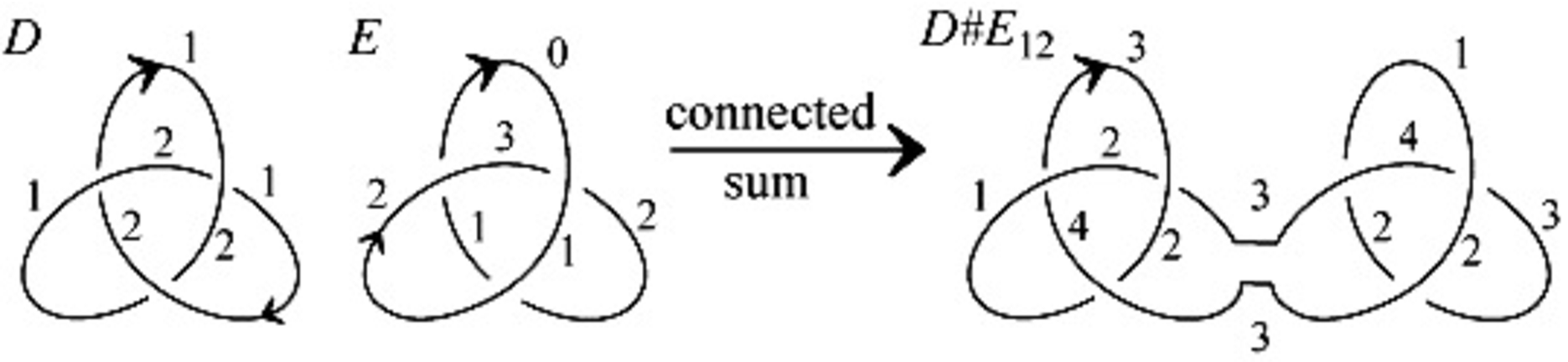}
\caption{}
\label{conn-ex}
\end{center}
\end{figure}

In this subsection, we discuss the spans of connected sums. 
Let $D$ and $E$ be knot diagrams on $S^2$ with warping degree labeling. 
Let $D\sharp E_{ij}$ be the connected sum of $D$ and $E$ at an edge of $D$ labeled $i$ and that of $E$ labeled $j$ 
as shown in Figure \ref{conn-ex}. 
We have the following proposition:

\phantom{x}
\begin{prop}
We have 
$$W_{D\sharp E_{ij}}(t)=t^jW_D(t)+t^iW_E(t).$$

\phantom{x}
\begin{proof}
An edge of $D\sharp E_{ij}$ corresponding to that of $D$ (resp. $E$) labeled $n$ has the warping degree $n+j$ (resp. $n+i$). 
\end{proof}
\label{conn-poly}
\end{prop}
\phantom{x}

\noindent We have the following Lemma: 

\phantom{x}
\begin{lem}
Let $D$ and $E$ be knot diagrams. We have 
$$\max (\mathrm{spn}(D), \mathrm{spn}(E))\le \mathrm{spn}(D\sharp E_{ij})\le \mathrm{spn}(D)+\mathrm{spn}(E).$$
Further, the first equality holds if and only if $i$ and $j$ satisfy 
$$\min _b(D_b)-\min _b(E_b)\le i-j\le \max _b(D_b)-\max _b(E_b),$$
where we assume that $\mathrm{spn}(D)\ge \mathrm{spn}(E)$. 

\phantom{x}
\begin{proof}
We have the first inequality and the inequality $\mathrm{spn}(D\sharp E_{ij})\le \mathrm{spn}(D)+\mathrm{spn}(E)+1$ by Proposition \ref{conn-poly}. 
Both $t^jW_D(t)$ and $t^iW_E(t)$ have terms with degree $i+j$ because $W_D(t)$ (resp. $W_E(t)$) has a term with degree $i$ (resp. $j$). 
Hence we have the second inequality. \\
Let $W_D(t)=n_0t^d+n_1t^{d+1}+\dots +n_st^{d+s}$ and $W_E(t)=m_0t^e+m_1t^{e+1}+\dots +m_rt^{e+r}$. Then, 
\begin{align*}
W_{D\sharp E_{ij}}(t)=&t^jW_D(t)+t^iW_E(t)\\
                          =&(n_0t^{d+j}+n_1t^{d+j+1}+\dots +n_st^{d+j+s})\\
                            &+(m_0t^{e+i}+m_1t^{e+i+1}+\dots +m_rt^{e+i+r}).
\end{align*}
Hence we have $\mathrm{spn}(D)=\mathrm{spn}(D\sharp E_{ij})$ if and only if $d+j\le e+i$ and $e+i+r\le d+j+s$, namely 
$d-e\le i-j\le (d+s)-(e+r)$. 
\end{proof}
\end{lem}
\phantom{x}

\begin{figure}[!ht]
\begin{center}
\includegraphics[width=110mm]{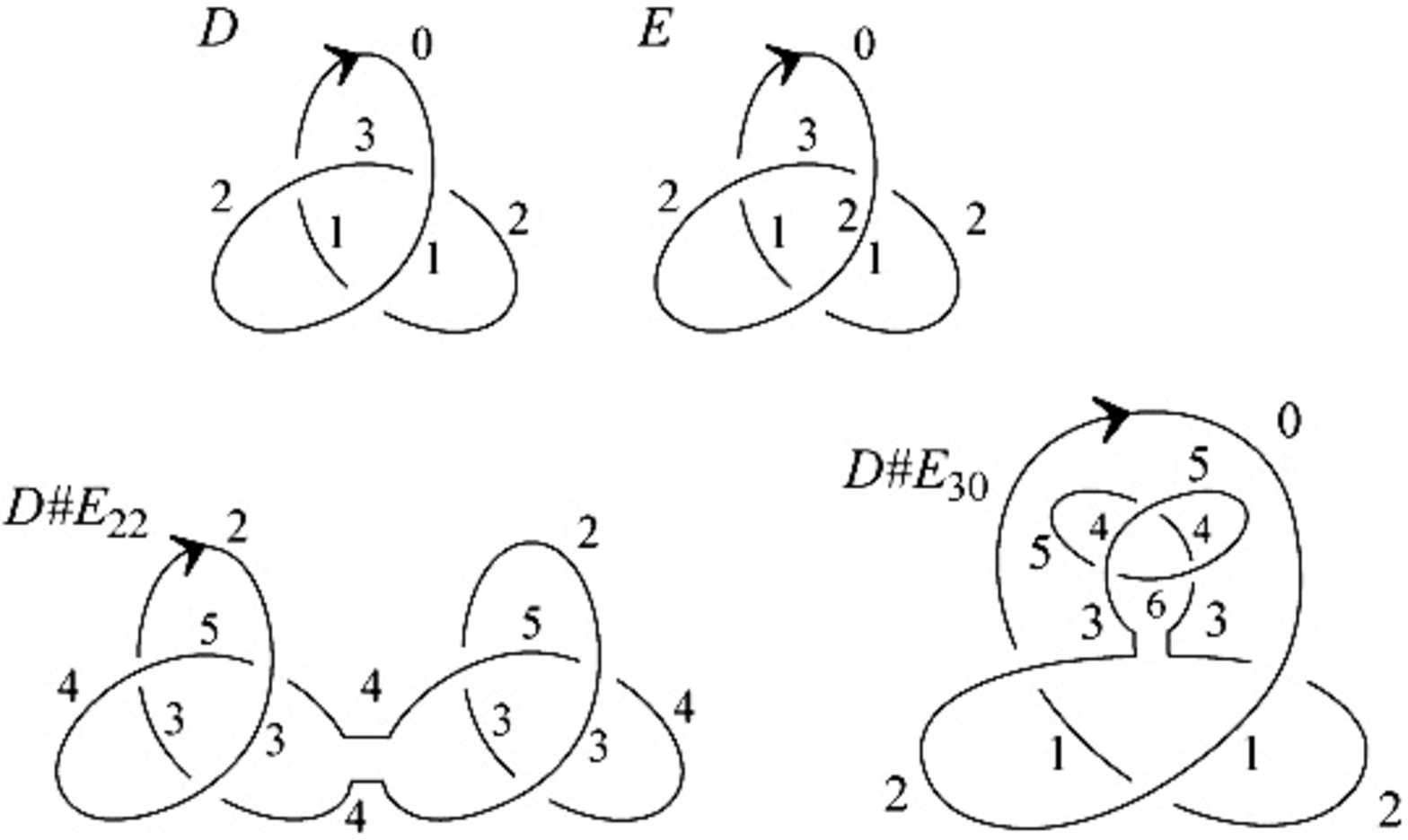}
\caption{}
\label{conn-ex2}
\end{center}
\end{figure}

\noindent For example, for two diagrams $D$, $E$ with $\mathrm{spn}(D)=\mathrm{spn}(E)=3$ in Figure \ref{conn-ex2}, we have 
$\mathrm{spn}(D\sharp E_{22})=3$ and $\mathrm{spn}(D\sharp E_{30})=6$. 
For knots, we have the following proposition: 

\phantom{x}
\begin{prop}
Let $K$ and $L$ be knots, and let $K\sharp L$ be the connected sum of $K$ and $L$. 
We have 
$$\mathrm{spn}(K\sharp L)\le \max (\mathrm{spn}(K), \mathrm{spn}(L)).$$

\phantom{x}
\begin{proof}
Let $D$, $E$ be diagrams of $K$, $L$ such that $\mathrm{spn}(D)=\mathrm{spn}(K)$, $\mathrm{spn}(E)=\mathrm{spn}(L)$. 
By choosing suitable $i$ and $j$, we have $\max (\mathrm{spn}(D), \mathrm{spn}(E))=\mathrm{spn}({D\sharp E}_{ij})\ge  \mathrm{spn}(K\sharp L)$. 
Hence $\max (\mathrm{spn}(K), \mathrm{spn}(L))\ge \mathrm{spn}(K\sharp L)$. 
\end{proof}
\end{prop}
\phantom{x}

\noindent We raise the following question:

\phantom{q}
\begin{q}
Let $K$ and $L$ be knots. Is it true that 
$$\mathrm{spn}(K\sharp L)= \max (\mathrm{spn}(K), \mathrm{spn}(L))?$$
\end{q}
\phantom{x}

\section*{Acknowledgments}
The author thanks Professor Akio Kawauchi for his valuable advice and encouregement. 
She also thanks Professors Taizo Kanenobu, Makoto Sakuma, Kengo Kishimoto, Yoshiro Yaguchi, and the members of Friday Seminar on Knot Theory 2010 at Osaka City University for valuable discussions and advice. 
She is deeply grateful to Professor Seiichi Kamada for his helpful suggestions including of the definition of the warping polynomial. 
She also thanks Professor Yasutaka Nakanishi for suggesting improvements of the previous version of this paper. 
She is partly supported by JSPS Research Fellowships for Young Scientists.

\maketitle

\end{document}